\documentclass[preprint,3p,12pt]{elsarticle}

\usepackage{lineno,hyperref}
\modulolinenumbers[5]

\usepackage{hyperref}

\usepackage{epstopdf}

\usepackage{mathrsfs}
\usepackage{amssymb}
\usepackage{amsfonts}
\usepackage{latexsym}
\usepackage{amsmath}
\usepackage{xcolor}
\usepackage{wrapfig}
\usepackage{floatflt}

\usepackage{mathtools}
\usepackage{extarrows}

\usepackage{graphicx}
\def\lb{\label}

\newcommand{\er}[1]{\textrm{(\ref{#1})}}

\newtheorem{theorem}{\bf Theorem}[section]

\def\a{\alpha}  \def\cA{{\mathcal A}}       
\def\b{\beta}   \def\cB{{\mathcal B}}       \def\mB{{\mathscr B}}
         
  \def\cD{{\mathcal D}}       
\def\d{\delta}         
  \def\cF{{\mathcal F}}       
           
          \def\mH{{\mathscr H}}

 \def\cM{{\mathcal M}}

    \def\cU{{\mathcal U}}       \def\mU{{\mathscr U}}
           \def\mV{{\mathscr V}}

       \def\vp{\varphi}    

\def\Z{{\mathbb Z}}    \def\R{{\mathbb R}}   \def\C{{\mathbb C}}    
\def\T{{\mathbb T}}    \def\N{{\mathbb N}}   

\def\lt{\biggl}                  \def\rt{\biggr}
\def\ol{\overline}

\def\iy{\infty}
               
\def\ss{\subset}                 \def\ts{\times}
                
                 \def\ev{\equiv}
        
\def\el2{\ell^{\,2}}             \def\1{1\!\!1}

\newcommand{\ca}{\begin{cases}}
\newcommand{\ac}{\end{cases}}
\newcommand{\ma}{\begin{pmatrix}}
\newcommand{\am}{\end{pmatrix}}
\def\eq{\begin{equation}}
\def\qe{\end{equation}}
\def\[{\begin{equation}}
\def\]{\end{equation}}

\begin{document}

\begin{frontmatter}



\title{Finite PDEs and finite ODEs are isomorphic}

\date{\today}

\author
{Anton A. Kutsenko}

\address{Jacobs University, 28759 Bremen, Germany; email: akucenko@gmail.com}

\begin{abstract}
The standard view is that PDEs are much more complex than ODEs, but, as will be shown below, for finite derivatives this  is not true. We consider the $C^*$-algebras ${\mathscr H}_{N,M}$ consisting of $N$-dimensional finite differential operators with $M\times M$-matrix-valued bounded periodic coefficients. We show that any ${\mathscr H}_{N,M}$ is $*$-isomorphic to
the universal uniformly hyperfinite algebra (UHF algebra) 
$
 \bigotimes_{n=1}^{\infty}\mathbb{C}^{n\times n}.
$
This is a complete characterization of the differential algebras. In particular, for different $N,M\in\mathbb{N}$ the algebras ${\mathscr H}_{N,M}$ are topologically and algebraically isomorphic to each other. In this sense, there is no difference between multidimensional matrix valued PDEs ${\mathscr H}_{N,M}$ and one-dimensional scalar ODEs ${\mathscr H}_{1,1}$. Roughly speaking, the multidimensional world can be emulated by the one-dimensional one.
\end{abstract}



\end{frontmatter}



{\section{Introduction}\lb{sec1}}

There is an obvious difference between linear continuous ordinary differential equations and partial differential systems, both with non-constant periodic coefficients.  In general, while ODEs is a part of PDEs formally, often books are written either about ODEs or PDEs, see, e.g., \cite{TP2012,E2010} (at least, the titles). In the continuous case, there are many reasons for this separation. For example, there is no full analogue of the Picard-Lindel\"of theorem even for linear PDEs.  Nevertheless, we show that if we replace continuous derivatives by their discrete analogues then both generate the same algebra, namely the universal UHF $C^*$-algebra. In this sense, topological and algebraic properties of algebras of finite ODEs and PDEs are identical. Algebras of discrete and continuous PDEs have numerous applications including a development of symbolic and numerical solvers of various differential equations, see, e.g., \cite{Ro1,GK01,Bu1,GRR}. 

Let us briefly describe another motivation of the paper related to the problems of non-linearity. It is well-known that non-linear stochastic ODEs lead to the linear Fokker-Planck PDEs describing the probability density function of the solution of non-linear stochastic ODE. There is also a more simple explanation why non-linear ODEs can be written as linear PDEs. Consider, possibly non-linear, equation $\dot{\bf x}={\bf f}(t,{\bf x})$, ${\bf x}|_{t=0}={\bf x}_0$. Let ${\bf x}(t)$ be the solution of this equation. Let $u(t,{\bf x})=\d({\bf x}-{\bf x}(t))$ be a formal ``trajectory'' of the solution in the phase-space, where $\d$ is a smooth approximation of the Kroneker delta. Then, formally differentiating $u$ we obtain the linear PDE
$$
 \dot u=-{\bf f}\cdot\nabla_{\bf x} u
$$
with the initial data $u(0,{\bf x})=\d({\bf x}-{\bf x}_0)$.
In this sense, the theory of non-linear ODEs is a part of the theory of linear PDEs. As mentioned above, we will show that the theory of finite linear PDEs is equivalent to the theory of finite linear ODEs. Roughly speaking, this means that finite analogues of non-linear ODEs, stochastic ODEs, and linear ODEs are more or less of the same type of complexity.

Let $N,M\in\N$ be positive integers. Let $L^2_{N,M}=L^2(\T^N\to\C^M)$ be the Hilbert space of periodic vector valued functions defined on the multidimensional torus $\T^N$, where $\T=\R/\Z\simeq[0,1)$. Everywhere in the article, it is assumed the Lebesgue measure in the definition of Hilbert spaces of square-integrable functions. Let $R^{\iy}_{N,M}=R^{\iy}(\T^N\to\C^{M\ts M})$ be the $C^*$-algebra of matrix-valued regulated functions with rational discontinuities. The regulated functions with possible rational discontinuities are the functions that can be uniformly approximated by the step functions of the form
\[\lb{REG}
 {\bf S}({\bf x})=\sum_{n=1}^P \chi_{J_n}({\bf x}){\bf S}_{n},
\]
where $P\in\N$, ${\bf S}_{n}\in\C^{M\ts M}$, and $\chi_{J_n}$ is the characteristic function of the parallelepiped $J_n=\prod_{i=1}^N[p_{in},q_{in})$ with rational end points $p_{in},q_{in}\in\mathbb{Q}/\Z\ss\T$. In particular, continuous matrix-valued functions belong to $R^{\iy}_{N,M}$. 

For ${\bf S}\in R^{\iy}_{N,M}$, the operator of multiplication by the function $\cM_{\bf S}: L^2_{N,M}\to L^2_{N,M}$ is defined by
\[\lb{multop}
 \cM_{\bf S}{\bf u}({\bf x})={\bf S}({\bf x}){\bf u}({\bf x}),\ \ {\bf u}({\bf x})\in L^2_{N,M},\ \ {\bf x}\in\T^N.
\]
For $i\in\N_N=\{1,...,N\}$, $h\in\mathbb{Q}$, the finite derivative $\cD_{i,h}:L^2_{N,M}\to L^2_{N,M}$ is defined by
\[\lb{diffop}
 \cD_{i,h}{\bf u}({\bf x})=\frac{{\bf u}({\bf x}+h{\bf e}_i)-{\bf u}({\bf x})}h,\ \ {\bf u}({\bf x})\in L^2_{N,M},\ \ {\bf x}\in\T^N,
\]
where the standard basis vector ${\bf e}_i=(\d_{ij})_{j=1}^N$, and $\d_{ij}$ is the Kronecker delta. The finite partial differential operators with bounded coefficients have the usual form
\[\lb{PDE}
 \cA {\bf u}=\sum_{n=1}^P\lt(\prod_{j=1}^{P_n}\cM_{jn}\cD_{jn}\rt){\bf u}+\cM_{00}{\bf u},\ \ {\bf u}\in L^2_{N,M},
\]
where $P,P_n\in\N$, and $\cM_{jn}$, $\cD_{jn}$ are some operators of the form \er{multop}, \er{diffop} respectively. 
The algebra of finite PDEs $\mH_{N,M}$ is generated by all the possible operators $\cA$ given by \er{PDE}, i.e.
\[\lb{APDE}
 \mH_{N,M}=\ol{{\rm Alg}}^{\mB}\{\cD_{i,h},\ \cM_{\bf S}:\ i\in\N_N,\ h\in\mathbb{Q},\ \ {\bf S}\in R^{\iy}_{N,M}\},
\] 
where $\mB=\mB_{N,M}=\mB(L^2_{N,M})$ is the $C^*$-algebra of bounded operators acting on $L^2_{N,M}$. It is seen that $\mH_{N,M}$ is the unital $C^*$-sub-algebra in $\mB$.

Let us introduce the universal UHF $C^*$-algebra $\mU$. One of the definitions is based on the inductive limit
$$
 \C^{1\ts1}\xrightarrow[]{\ \vp_1\ }\C^{2\ts2}\xrightarrow[]{\ \vp_2\ }\C^{6\ts6}...\xrightarrow[]{\ \vp_{n-1}\ }\C^{n!\ts n!}...\xrightarrow[]{}\mU,
$$
where $\vp_n$ are the unital $*$-embeddings, or, formally,
$$
 \mU=\C^{1\ts1}\otimes\C^{2\ts2}\otimes\C^{3\ts3}\otimes...\C^{n\ts n}\otimes....
$$ 
The corresponding supernatural number $\mathfrak{N}_{\mU}$ contains all prime numbers infinitely many times. Hence, any UHF algebra is a sub-algebra of $\mU$, since any supernatural number devides $\mathfrak{N}_{\mU}$. Recall that there is one-to-one correspondence between UHF algebras and supernatural numbers, see, e.g., \cite{G1960,D1997,RLL2000}. Let us formulate our main result.

\begin{theorem}\lb{T1} For any $N,M\in\N$, the $C^*$-algebra $\mH_{N,M}$ is $*$-isomorphic to $\mU$. Moreover, there is a unitary $\cU_{N,M}:L^2_{1,1}\to L^2_{N,M}$ such that $\mH_{1,1}=\cU^{-1}_{N,M}\mH_{N,M}\cU_{N,M}$.
\end{theorem}
 
This means that there is no difference between $\mH_{N,M}$ for different $N,M$. For example, if $\cA\in\mH_{N,M}$ then there is $\cB\in\mH_{1,1}$ with the same spectrum ${\rm sp}_{\mB_{N,M}}(\cA)={\rm sp}_{\mB_{1,1}}(\cB)$,
and the $C^{*}$-algebras generated by $\cA$ and $\cB$ are $*$-isomorphic
$$
 \ol{{\rm Alg}}^{\mB_{N,M}}\{1,\cA,\cA^*\}\cong\ol{{\rm Alg}}^{\mB_{1,1}}\{1,\cB,\cB^*\}.
$$ 
In particular, $\cA$ is invertible if and only if $\cB$ is invertible.  Thus, there are no special difficulties in the analysis of finite PDEs in comparison with finite ODEs. Moreover, there is a unitary transform between solutions of ODEs and PDEs given by the unitary operator $\cU$, see  Theorem \ref{T1}.

It is useful to take into account the following remark. Using  $L^2_{N,M}=\bigoplus_{m=1}^M(L^2_{1,1})^{\otimes N}$, we conclude that
\[\lb{1dtoNd}
 \mH_{N,M}\cong\C^{M\ts M}\otimes\mH_{1,1}^{\otimes N}\cong\C^{M\ts M}\otimes\mU^{\otimes N}\cong\mU\cong\mH_{1,1},
\]
since $\mU\otimes\mV\cong\mU\cong\mV\otimes\mU$ for any UHF algebra $\mV$. 

What about other UHF algebras? Let $\mathfrak{N}_{\mV}=\prod_{n=1}^{P}p_n^{N_n}$ be some supernatural number corresponding to the UHF algebra $\mV$. Some of $N_n$ and $P$ can be infinite. In addition to the notation $\N_N=\{1,...,N\}$, we use also $\N_{\iy}=\N$. Then
\[\lb{ISO11}
 \mV\cong\mH_{1,1}(\mathfrak{N}_{\mV})\ev\ol{{\rm Alg}}^{\mB_{1,1}}\{\cD_{1,h},\ \cM_{ S}:\ h=p_n^{-j},\ j\in\N_{N_n},\ n\in\N_P,\ { S}\in R^{\iy}_{1,1}(\mathfrak{N}_{\mV})\},
\]
where $R^{\iy}_{1,1}(\mathfrak{N}_{\mV})$ is the $C^{*}$-algebra of scalar one dimensional regulated functions which can be uniformly approximated by step functions of the form \er{REG} but with $q_{1n},p_{1n}$ equal to some $p_r^{-s}$ for any $r\in\N_P$, $s\in\N_{N_r}$.  In particular, the CAR-algebra (canonical anticommutation relations in quantum mechanics) admits the representation as the differential algebra generated by the dyadic derivatives $\cD_{i,2^{-n}}$, $n\in\N$ and dyadic regulated functions. The proof of \er{ISO11} is similar to the proof of Theorem \ref{T1}. 

Let $\mathfrak{N}$ be some supernatural number.  For the multidimensional case, we define
\[\lb{ISONM}
 \mH_{N,M}(\mathfrak{N})=\C^{M\ts M}\otimes\mH_{1,1}(\mathfrak{N})^{\otimes N},
\]
see the first identity in \er{1dtoNd}. Then the corresponding supernatural number is
$$
 \mathfrak{N}_{\mH_{N,M}(\mathfrak{N})}=M\mathfrak{N}^N,
$$
since $\mathfrak{N}_{\mU\otimes\mV}=\mathfrak{N}_{\mU}\mathfrak{N}_{\mV}$ for any UHF algebras $\mU,\mV$. For example, $\mH_{N,M}(\mathfrak{N})$ is the CAR-algebra if and only if $\mathfrak{N}=2^{\iy}$ and $M=2^m$ for some $m\in\N\cup\{0\}$.

{\section{Proof of Theorem \ref{T1}}\lb{sec2}}

We fix $N,M\in\N$ and, for convenience, we will omit these indices below. Let $h=1/p$ for some $p\in\N$. Denote $R^{h}$ the $C^*$-sub-algebra of $R^{\iy}$ consisting of step functions constant on each $J=\prod_{i=1}^N[hp_i,hp_i+h)\ss\T^N$, where $p_i\in\Z_p=\Z/p\Z=\{0,...,p-1\}$. Consider the finite-dimensional $C^*$-sub-algebra $\mH^h\ss\mH$ defined by
$$
 \mH^h={\rm Alg}\{\cD_{i,h},\ \cM_{\bf S}:\ i\in\N_N,\ \ {\bf S}\in R^{h}\}.
$$
It is seen that any operator $\cA\in\mH^h$ has the form
\[\lb{102}
 \cA{\bf u}({\bf x})=\sum_{{\bf j}\in\Z_p^N}{\bf A}_{{\bf j}}({\bf x}){\bf u}({\bf x}+h{\bf j}),\ \ {\bf u}\in L^2_{N,M},
\]
with some ${\bf A}_{{\bf j}}({\bf x})\in R^h$. This is because all such $\cA$ belongs to $\mH^h$, since $\cD_{i,h}$ is generated by shift operators, and all such $\cA$ form an algebra which contains $\cD_{i,h}$ and $\cM_{\bf S}$ for ${\bf S}\in R^h$.
The next step is to find the convenient representation of $\mH^h$. Using \er{102}, we have
\begin{multline}\lb{103}
 \cA{\bf u}({\bf y}+h{\bf r})=\sum_{{\bf j}\in\Z_n^N}{\bf A}_{{\bf j}}({\bf y}+h{\bf r}){\bf u}({\bf y}+h{\bf j}+h{\bf r})=
 \sum_{{\bf j}\in\Z_n^N}{\bf A}_{{\bf j}-{\bf r}}({\bf y}+h{\bf r}){\bf u}({\bf y}+h{\bf j}),\\ {\bf y}\in I_h=[0,h)^N,\ \ {\bf r}\in\Z_p^N.
\end{multline}
The Hilbert space $L^2$ is naturally isomorphic to the direct sum of Hilbert spaces of functions defined on shifted cubes $I_h$:
$$
 L^2\cong L_h^2=\bigoplus_{{\bf r}\in\Z_p^N}L^2(I_h\to\C^M)
$$
with the isomorphism $\cF_h: L^2\to L^2_h$ defined by
\[\lb{104}
 {\bf u}({\bf x})\longleftrightarrow({\bf u}({\bf y}+h{\bf j}))_{{\bf j}\in\Z_p^N},\ \ {\bf u}\in L^2,\ \ {\bf x}\in\T^N,\ \ {\bf y}\in I_h.
\]
Then, by \er{103} the operator $\cF_h\cA\cF_h^{-1}:L^2_h\to L^2_h$ is the operator of multiplication by $Mp^N\ts Mp^N$-matrix-valued function ${\bf B}_{\cA}({\bf y})$, ${\bf y}\in I_h$ of the form
\[\lb{105}
 {\bf B}_{\cA}({\bf y})=({\bf A}_{{\bf j}-{\bf r}}({\bf y}+h{\bf r}))_{{\bf r},{\bf j}\in\Z_p^N}.
\]
It is easy to see that
\[\lb{106}
 {\bf B}_{\a\cA+\b\cB}=\a{\bf B}_{\cA}+\b{\bf B}_{\cB},\ \ {\bf B}_{\cA\cB}={\bf B}_{\cA}{\bf B}_{\cB},\ \ {\bf B}_{\cA^*}={\bf B}^*_{\cA}.
\]
Moreover, ${\bf B}_{\cA}({\bf y})$ is constant for ${\bf y}\in I_h$, since all ${\bf A}_{\bf j}\in R^h$. Finally, note that for any ${\bf B}\in\C^{Mp^N\ts Mp^N}$ there is a unique $\cA\in\mH^h$, such that ${\bf B}={\bf B}_{\cA}$. We can explicitly and uniquely recover $\cA$ from ${\bf B}_{\cA}$, using \er{105} and \er{103}. Thus, $\cA\mapsto{\bf B}_{\cA}$ is $*$-isomorphism between $\mH^h$ and $\C^{Mp^N\ts Mp^N}$.

Taking $h_n=1/n!$, $n\in\N$, we can write
\[\lb{107}
 \mH^{h_1}\xrightarrow[]{\ \vp_1\ }\mH^{h_2}\xrightarrow[]{\ \vp_3\ }\mH^{h_3}\xrightarrow[]{\ \vp_2\ }...\xrightarrow[]{}\mH,
\]
where $\vp_n$ is the natural embedding $\mH^{h_n}$ into $\mH^{h_{n+1}}$. Such embedding exists, since the partition of $\T^N$ onto $((n+1)!)^N$ identical cubes for $\mH^{h_{n+1}}$ contains the partition of $\T^N$ onto $(n!)^N$ identical cubes for $\mH^{h_{n}}$. The inductive limit in \er{107} is $\mH$, since any $\cD_{i,h}$ with $h=p/q\in\mathbb{Q}$ ($q\in\N$) belongs to $\mH^{h_q}$, and $R^{\iy}$ can be uniformly approximated by $R^{h_n}$ following the definition of regulated functions. Remembering $\mH^{h_n}\cong\C^{M(n!)^N\ts M(n!)^N}$, we can conclude that the supernatural number $\mathfrak{N}_{\mH}$ for the UHF algebra $\mH$ contains all prime numbers infinitely many times. Hence, $\mH$ is $*$-isomorphic to the universal UHF algebra $\mU$.

\begin{figure}[h]
	\center{\includegraphics[width=0.65\linewidth]{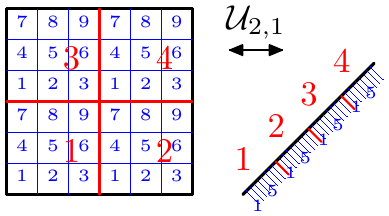}}
	\caption{Two first partitions of the unitary transform $\cU_{2,1}$ between $L^2_{2,1}$ and $L^2_{1,1}$ are shown. The characteristic functions of squares and intervals with the same "blue" and "red" numbers are transformed into each other.}
	\label{fig1}
\end{figure}

Let us construct the unitary operator $\cU:L^2_{1,1}\to L^2_{N,M}$ such that $\mH_{1,1}=\cU^{-1}\mH_{N,M}\cU$. The unitary $\cU$ can be any operator which transform characteristic functions of cubic cells to characteristic functions of intervals preserving the order, see Fig. \ref{fig1}. Any number $x\in[0,1)$ can be expanded as
\[\lb{exp1}
 x=\frac{x_1}{M}+\frac{x_2}{M(2!)^N}+\frac{x_3}{M(3!)^N}+...,
\]
where
\[\lb{exp2} 
 x_1\in\{0,...,M-1\},\ \ \ x_i\in\{0,...,i^N-1\},\ i>1.
\] 
The coefficients $x_i$ can be found recurrently
\[\lb{exp3}
 x_1=\lfloor Mx\rfloor,\ \ x_2=\lfloor (Mx-x_1)2^N\rfloor,\ \ x_3=\lfloor (M2^Nx-2^Nx_1-x_2)3^N \rfloor,\ \ ....
\] 
Now, let
\[\lb{exp4}
 \vp_i:\{0,...,i^N-1\}\to\{0,...,i-1\}^N
\]
be some 1-1 mappings. Define the mapping $\pmb{\vp}:[0,1/M)\to[0,1]^N$ by
\[\lb{exp5}
 \pmb{\vp}(x)=\frac{\vp_2(x_2)}{2!}+\frac{\vp_3(x_3)}{3!}+\frac{\vp_4(x_4)}{4!}+....
\]
Then, we can define the unitary $\cU$ by
\[\lb{exp6}
 \cU u={\bf v}=(v_{i+1})_{i=0}^{M-1},\ \ where\ \ u(x)|_{[\frac iM,\frac{i+1}M)}=M^{-\frac12}v_{i+1}(\pmb{\vp}(x-\frac iM))
\]
and $u\in L^2_{1,1}$, ${\bf v}\in L^2_{N,M}$. We need the factor $M^{-\frac12}$ because $\cU$ should be unitary operator. Note also that \er{exp6} is valid for almost all $x$ except some set of zero measure. It's because of the fact that $\pmb{\vp}$ is an injection up to a set of zero measure. This is an analog of the fact that  digital expansion is unique to all the real numbers except some set of zero measure.

\ 

\section*{Acknowledgements} This work is supported by the Russian Science Foundation (RSF) project 18-11-00032. This paper is also a contribution to the project M3 of the Collaborative Research Centre TRR 181 ``Energy Transfer in Atmosphere and Ocean'' funded by the Deutsche Forschungsgemeinschaft
(DFG, German Research Foundation) under project number 274762653.

\medskip

\bibliography{bibl_perp1}

\end{document}